\documentclass[a4paper, 12pt]{amsart}
\usepackage{amsmath,amssymb,amsfonts}
\usepackage{mathrsfs,latexsym,amsthm,stmaryrd,enumerate}
\usepackage[all]{xy}
\usepackage[vcentermath]{youngtab}

\newcommand{\bbC}{\mathbb{C}}

\newcommand{\tto}{\twoheadrightarrow}
\renewcommand{\to}{\rightarrow}

\renewcommand{\phi}{\varphi}
\renewcommand{\epsilon}{\varepsilon}

\newtheorem{proposition}{Proposition}
\newtheorem{lemma}[proposition]{Lemma}
\newtheorem{theorem}[proposition]{Theorem}
\newtheorem{corollary}[proposition]{Corollary}
\newtheorem{conjecture}[proposition]{Conjecture}

\theoremstyle{definition}
\newtheorem{remark}[proposition]{Remark}

\def\clap#1{\hbox to 0pt{\hss#1\hss}}

\begin{document}

\title{A new approach to Kostant's problem}
\author{Johan K{\aa}hrstr{\"o}m and Volodymyr Mazorchuk}
\date{\today}
\maketitle

\begin{abstract}
For every involution $\mathbf{w}$ of the symmetric group
$S_n$ we establish, in terms ofa special canonical quotient of the 
dominant Verma module associated with $\mathbf{w}$, an effective
criterion, which allows us to verify whether the universal enveloping 
algebra  $U(\mathfrak{sl}_n)$ surjects onto the space of all ad-finite 
linear transformations of the simple highest weight module $L(\mathbf{w})$.
An easy sufficient condition derived from this criterion admits 
a straightforward computational check for example using a computer.
All this is applied to get some old and many new results,
which answer the classical question of Kostant in special cases,
in particular we give a complete answer for simple highest weight
modules in the regular block of $\mathfrak{sl}_n$, $n\leq 5$.  
\end{abstract}

\section{Introduction}\label{s0}

Let $\mathfrak{g}$ be a complex semi-simple finite-dimensional Lie algebra 
with a fixed triangular decomposition, 
$\mathfrak{g}=\mathfrak{n}_-\oplus \mathfrak{h}\oplus\mathfrak{n}_+$, and 
$U(\mathfrak{g})$ be its universal enveloping algebra. Then for every
two $\mathfrak{g}$-modules $M$ and $N$ the space 
$\mathrm{Hom}_{\mathbb{C}}(M,N)$ may be viewed as a 
$U(\mathfrak{g})$-bimodule in the natural way, and, furthermore, 
also as a $\mathfrak{g}$-module under the adjoint action of $\mathfrak{g}$. 
The bimodule $\mathrm{Hom}_{\mathbb{C}}(M,N)$ has a sub-bimodule, usually
denoted by $\mathcal{L}(M,N)$, which consists of all elements, on which 
the adjoint action of $U(\mathfrak{g})$ is locally finite  (see for example \cite[Kapitel~6]{Ja}). Since $U(\mathfrak{g})$ itself consists of locally 
finite elements under the adjoint action, it naturally maps to 
$\mathcal{L}(M,M)$ for every $\mathfrak{g}$-module $M$, and the kernel of 
this map is obviously the annihilator $\mathrm{Ann}(M)$ of $M$ in
$U(\mathfrak{g})$. The classical problem of Kostant (see 
for example \cite{Jo2}) is formulated  in the following way: 
\vspace{4mm}

\noindent
{\em For which $\mathfrak{g}$-modules $M$ is the natural injection 
\begin{displaymath}
U(\mathfrak{g})/\mathrm{Ann}(M)\hookrightarrow \mathcal{L}(M,M)
\end{displaymath}
surjective?}
\vspace{4mm}

The (positive) answer to Kostant's problem is an important tool, 
in particular, in the study of generalized Verma modules, 
see \cite{MiSo,KM1,MS1}. Unfortunately, the complete answer to this 
problem is not even known for simple highest weight modules. 
The answer is known to be positive for Verma 
modules (see \cite[Corollary~6.4]{Jo2}) and for certain classes of 
simple highest weight modules (see \cite[Theorem~4.4]{GJ} and
\cite[Theorem~1]{Ma}). For simple highest weight modules in type $A$ the 
answer is even known to be an invariant of a left cell, see 
\cite[Theorem~60]{MS1}. However, already in \cite[9.5]{Jo2} it was shown 
that for some simple highest weight modules in type $B$ the answer is 
negative. In spite of the general belief that the answer is positive
for simple highest weight modules in type $A$, it was recently shown in
\cite[Theorem~13]{MS2} that for the simple highest weight 
$\mathfrak{sl}_4$-module $L(rt)$, where $r$ and $t$ are two commuting
simple reflections, the answer is negative.

The present paper is strongly inspired by the latter counter-example
and is an attempt to analyze and generalize it. As in type $A$ the 
answer to Kostant's problem is an invariant of a left cell, and since
every left cell of the symmetric group $S_n$ contains a unique
involution, it is enough to solve Kostant's problem for all modules 
of the form $L(\mathbf{w})$, where  $\mathbf{w}\in S_n$ is an involution. 
The counter-example in \cite[Theorem~13]{MS2} was constructed relating the
module $L(\mathbf{w})$ to a special quotient of the dominant Verma module,
which in the following will be denoted by $D^{\hat{\mathrm{R}}}$. This
module is a canonical object of the category 
$\mathcal{O}_0^{\hat{\mathrm{R}}}$, which was used in \cite{MS1}
to categorify Kazhdan-Lusztig cell modules. The module $L(\mathbf{w})$
is the simple socle of $D^{\hat{\mathrm{R}}}$ and thus both 
$L(\mathbf{w})$ and $D^{\hat{\mathrm{R}}}$ are submodules of the
indecomposable injective  module $P^{\hat{\mathrm{R}}}(\mathbf{w})$
in $\mathcal{O}_0^{\hat{\mathrm{R}}}$, which also turns out to be projective. 

The main result of the present paper relates the solution of Kostant's
problem for $L(\mathbf{w})$ to the structure of $D^{\hat{\mathrm{R}}}$
as follows: 

\begin{theorem}\label{mainresult}
Kostant's problem has a positive answer for $L(\mathbf{w})$ if and only if
every simple submodule of the cokernel of the canonical inclusion
$D^{\hat{\mathrm{R}}}\subset P^{\hat{\mathrm{R}}}(\mathbf{w})$ has the
form $L(x)$, where $x$ is some element from the right cell of $\mathbf{w}$.
\end{theorem}

We will show that Theorem~\ref{mainresult} can be used to answer
Kostant's problem in many cases, in particular, to obtain many new
results and reprove some old results. The most interesting application
of this theorem seems to be that it implies a sufficient condition
for a {\em negative} answer to Kostant's problem, which is purely
computational and can be realized as a relatively short and efficient
program on a computer.

In Section~\ref{s1} we collected all necessary preliminaries. The
main results (in particular Theorem~\ref{mainresult}) are formulated
in detail and proved in Section~\ref{s2}. In Section~\ref{s3} we
collected many applications, both theoretical and computational.
\vspace{1cm}

\noindent
{\bf Acknowledgments.} For the second author the research was
partially supported by the Swedish Research Council.

\section{Notation and preliminaries}\label{s1}

From now on we assume that $\mathfrak{g}=\mathfrak{sl}_n$ and the triangular
decomposition is just the usual decomposition into the upper triangular,
diagonal and lower triangular matrices. The symmetric group $S_n$ is the
Weyl group $W$ for $\mathfrak{g}$ and hence $S_n$ acts on $\mathfrak{h}^*$
in the usual way $w\lambda$, and via the dot action 
$w\cdot \lambda=w(\lambda+\rho)-\rho$, where $\rho$ is half the sum of all
positive (with respect to the above triangular decomposition) roots of the algebra $\mathfrak{g}$. 

Let $\mathcal{O}$ denote the BGG category $\mathcal{O}$, \cite{BGG},
associated with the triangular decomposition above. For $w\in W$ we let
$\Delta(w)$ denote the Verma module with highest weight $w\cdot 0$, $L(w)$
denote the simple head of $\Delta(w)$, and $P(w)$ denote the indecomposable
projective cover of $L(w)$. The {\em principal block} $\mathcal{O}_0$ of
$\mathcal{O}$ is the full subcategory of $\mathcal{O}$, which contains all
$L(w)$, $w\in S_n$, and is closed under isomorphisms and extensions. 
The category $\mathcal{O}_0$ is a direct summand of $\mathcal{O}$.

For $w\in W$ we denote by $\theta_w$ the indecomposable {\em projective
functor} on $\mathcal{O}_0$, associated with $w$. This functor is a
unique (up to isomorphism) indecomposable direct summand of all possible 
functors, which have the form 
$V\otimes_{\mathbb{C}}{}_-:\mathcal{O}\to \mathcal{O}$,
where $V$ is a finite-dimensional $\mathfrak{g}$-module, which 
satisfies $\theta_w\Delta(e)=P(w)$, see \cite[Section~3]{BG}.

Denote by $\leq_{\mathtt{L}}$ and $\leq_{\mathtt{R}}$ the left and the right
(pre)orders on $W$ respectively, see \cite[Section~3]{BB}. For a fixed 
right cell $\mathbf{R}$ set
\begin{displaymath}
\hat{\mathbf{R}}=\{x\in W\,:\, x\leq_{\mathtt{R}} w\text{ for some }w\in
\mathbf{R} \} 
\end{displaymath}
and denote by $\mathcal{O}_0^{\hat{\mathbf{R}}}$ the full subcategory 
of $\mathcal{O}_0$, which contains all $L(w)$, $w\in \hat{\mathbf{R}}$, 
and is closed under isomorphisms and extensions. The natural inclusion
functor $\mathfrak{i}_{\hat{\mathbf{R}}}^0:\mathcal{O}_0^{\hat{\mathbf{R}}}
\to \mathcal{O}_0$ is obviously exact and hence has both the left adjoint
$\mathrm{Z}^{\hat{\mathbf{R}}}_0:\mathcal{O}_0\to 
\mathcal{O}_0^{\hat{\mathbf{R}}}$ and the right adjoint
$\hat{\mathrm{Z}}^{\hat{\mathbf{R}}}_0:\mathcal{O}_0\to 
\mathcal{O}_0^{\hat{\mathbf{R}}}$, see \cite[5.1]{MS1}. The functor
$\mathrm{Z}^{\hat{\mathbf{R}}}_0$ is just the functor of taking the maximal
possible quotient, which lies in $\mathcal{O}_0^{\hat{\mathbf{R}}}$; and 
the functor $\hat{\mathrm{Z}}^{\hat{\mathbf{R}}}_0$ is just the functor 
of taking the maximal possible submodule, which lies in 
$\mathcal{O}_0^{\hat{\mathbf{R}}}$. All projective functors on
$\mathcal{O}_0$ preserve $\mathcal{O}_0^{\hat{\mathbf{R}}}$, and both
$\mathrm{Z}^{\hat{\mathbf{R}}}_0$ and $\hat{\mathrm{Z}}^{\hat{\mathbf{R}}}_0$
commute with $\theta_w$ for all $w\in W$, see \cite[Lemma~19]{MS1}.

For $w\in \hat{\mathbf{R}}$ set $P^{\hat{\mathbf{R}}}(w)=
\mathrm{Z}^{\hat{\mathbf{R}}}_0 P(w)$ and $\Delta^{\hat{\mathbf{R}}}(w)=
\mathrm{Z}^{\hat{\mathbf{R}}}_0 \Delta(w)$. Then the modules
$P^{\hat{\mathbf{R}}}(w)$,  $w\in \hat{\mathbf{R}}$, are exactly the 
indecomposable projective modules in $\mathcal{O}_0^{\hat{\mathbf{R}}}$.
The module $P^{\hat{\mathbf{R}}}(w)$ is injective if and only if
$w\in \mathbf{R}$, see \cite[Section~5]{MS1}. Let $\mathbf{w}\in \mathbf{R}$
be a unique involution in $\mathbf{R}$. Then 
$P^{\hat{\mathbf{R}}}(w)=\theta_{w}L(\mathbf{w})$ for any $w\in \mathbf{R}$,
see \cite[Key statement]{MS2}. By \cite[Lemma~8]{MS2} we have the equality 
$\dim\mathrm{Hom}_{\mathfrak{g}}(P^{\hat{\mathbf{R}}}(e),
P^{\hat{\mathbf{R}}}(\mathbf{w}))=1$. Denote by $D^{\hat{\mathbf{R}}}$
the image of the unique (up to a scalar) non-zero homomorphism from
$P^{\hat{\mathbf{R}}}(e)$ to $P^{\hat{\mathbf{R}}}(\mathbf{w})$.

\begin{conjecture}\label{conj1}
$D^{\hat{\mathbf{R}}}=P^{\hat{\mathbf{R}}}(e)$ .
\end{conjecture}

Define the following full subcategories in $\mathcal{O}_0^{\hat{\mathbf{R}}}$:
\begin{gather*}
\mathcal{C}_1=\{M\in \mathcal{O}_0^{\hat{\mathbf{R}}}\,:\,
[M:L(x)]>0\text{ implies }x<_{\mathtt{R}}\mathbf{w}\},\\
\mathcal{C}_2=\{M\in \mathcal{O}_0^{\hat{\mathbf{R}}}\,:\,
\mathrm{Hom}_{\mathfrak{g}}(L(x),M)\neq 0\text{ implies }x\in\mathbf{R}\},\\
\mathcal{C}_3=\{M\in \mathcal{O}_0^{\hat{\mathbf{R}}}\,:\,
\mathrm{Hom}_{\mathfrak{g}}(M,L(x))\neq 0\text{ implies }x\in\mathbf{R}\}
\end{gather*}
From the definition we immediately have
$\mathrm{Hom}_{\mathfrak{g}}(M,N)=0$ for all $X\in \mathcal{C}_1$ and
$Y\in \mathcal{C}_2$; and for all $X\in \mathcal{C}_3$ and
$Y\in \mathcal{C}_1$. 

\begin{lemma}\label{lemintro}
For every $w\in W$ and $i=1,2,3$ the functor $\theta_w$ preserves
the category $\mathcal{C}_i$.
\end{lemma}

\begin{proof}
That $\theta_w$ preserves $\mathcal{C}_1$ follows from the definitions
and the fact that $\theta_w$ preserves $\mathcal{O}_0^{\hat{\mathbf{R}}}$.
That $\theta_w$ preserves $\mathcal{C}_2$ follows from the fact that 
the injective envelope of any $X\in \mathcal{C}_2$ is projective and the
fact that $\theta_w$ is exact and preserves projective-injective modules
in $\mathcal{O}_0^{\hat{\mathbf{R}}}$. The proof of the fact that 
$\theta_w$ preserves $\mathcal{C}_3$ is dual. 
\end{proof}

Let $\mathcal{P}=\oplus_{w\in \mathbf{R}}P^{\hat{\mathbf{R}}}(w)$.
For every $M\in \mathcal{O}_0^{\hat{\mathbf{R}}}$ let $I_M$ be some
injective envelope of $M$ and set
\begin{displaymath}
M_1=\bigcap_{\substack{f\in\mathrm{Hom}_{\mathfrak{g}}
(I_M,\mathcal{P})\\f(M)=0}}\mathrm{Ker}(f),\quad\quad\quad
M'_1=\bigcap_{f\in\mathrm{Hom}_{\mathfrak{g}}
(M_1,\mathcal{P})}\mathrm{Ker}(f), 
\end{displaymath}
and $M_2=M_1/M'_1$. The correspondence $M\mapsto M_2$ is functorial and
$M_2$ is called the {\em partial approximation} of $M$ with respect to 
the injective module $\mathcal{P}$, see \cite[2.4]{KM}. We denote by
$\mathrm{A}:\mathcal{O}_0^{\hat{\mathbf{R}}}\to
\mathcal{O}_0^{\hat{\mathbf{R}}}$ the corresponding functor of partial
approximation. From the definition we have the natural transformation
$\mathrm{nat}$ from the identity functor $\mathrm{ID}$ to $\mathrm{A}$,
which is just the quotient map from $M$ to $M/(M\cap M'_1)$. The functor
$\mathrm{A}$ is left exact, see \cite[2.4]{KM}.

\section{The main results}\label{s2}

\subsection{A criterion for testing Kostant's problem}\label{s2.1}

According to \cite[Theorem~60]{MS1}, the answer to Kostant's problem for
$L(w)$, $w\in W$, is an invariant of a left cell. Since every left cell
has a unique involution, it is thus enough to study Kostant's problem
for involutions in $W$. The main result of the paper is the following 
statement:

\begin{theorem}\label{thm2}
Let $\mathbf{w}\in W$ be an involution and $\mathbf{R}$ be the right
cell of $W$, containing $\mathbf{w}$. Then the following conditions 
are equivalent:
\begin{enumerate}[(a)]
\item\label{thm2.1} Kostant's problem has a positive solution for
$L(\mathbf{w})$. 
\item\label{thm2.2} Every simple module, occurring in the socle of the
cokernel $\mathrm{Coker}$ of the natural 
inclusion $D^{\hat{\mathbf{R}}}\hookrightarrow
P^{\hat{\mathbf{R}}}(\mathbf{w})$, has the form $L(x)$, where $x\in 
\mathbf{R}$ (i.e. $\mathrm{Coker}$ belongs to $\mathcal{C}_2$).
\end{enumerate}
\end{theorem}

The idea of the proof is to compare Kostant's problem for modules
$L(\mathbf{w})$ and $D^{\hat{\mathbf{R}}}$. The former is exactly the 
module for which we would like to solve Kostant's problem, while the
latter is, by definition, a quotient of $\Delta(e)$, and hence
Kostant's problem for it has a positive solution by \cite[6.9(10)]{Ja}.
The relation between these two modules is again given by definition:
$L(\mathbf{w})$ is the simple socle of $D^{\hat{\mathbf{R}}}$. So, to 
compare $\mathcal{L}(L(\mathbf{w}),L(\mathbf{w}))$ and
$\mathcal{L}(D^{\hat{\mathbf{R}}},D^{\hat{\mathbf{R}}})$ one might 
first try to show that these two modules have the same annihilators, 
and then try to show that 
\begin{equation}\label{eq3}
\mathrm{Hom}_{\mathfrak{g}}(L(\mathbf{w}),\theta_w L(\mathbf{w}))=
\mathrm{Hom}_{\mathfrak{g}}(D^{\hat{\mathbf{R}}},
\theta_w D^{\hat{\mathbf{R}}})
\end{equation}
for all $w\in W$. This would be enough to conclude that 
$\mathcal{L}(L(\mathbf{w}),L(\mathbf{w}))=
\mathcal{L}(D^{\hat{\mathbf{R}}},D^{\hat{\mathbf{R}}})$ by 
\cite[6.8(3)]{Ja}, thus solving positively Kostant's problem for 
$L(\mathbf{w})$. The best way to prove \eqref{eq3} would be to construct
a functor, which commutes with all $\theta_w$, and sends $L(\mathbf{w})$
to $D^{\hat{\mathbf{R}}}$. It turns out that the functor $\mathrm{A}$
defined above does this job. So now let's do the work.

\begin{lemma}\label{lem4}
For all $w\in W$ there is an isomorphism of functors as follows:
$\mathrm{A}\theta_w\cong \theta_w\mathrm{A}$.
\end{lemma}

\begin{proof}
As $\mathrm{A}$ is left exact and $\theta_w$ is exact, both 
$\mathrm{A}\theta_w$ and $\theta_w\mathrm{A}$ are left exact.

Let $I\in \mathcal{O}_0^{\hat{\mathbf{R}}}$ be injective. Consider the
short exact sequence
\begin{equation}\label{eq5}
0\to K\to I\overset{\mathrm{nat}_I}{\longrightarrow} \mathrm{A}I\to 0,
\end{equation}
where $K$ is just the kernel of $\mathrm{nat}_I$. Since the socle of
$\mathcal{P}$ coincides with $\oplus_{w\in \mathbf{R}}L(w)$, from the
definition of $\mathrm{A}$ we have that $K\in\mathcal{C}_1$, while
$\mathrm{A}I\in\mathcal{C}_2$.

Applying $\theta_w$ to \eqref{eq5} and using Lemma~\ref{lemintro}
we obtain that $\theta_w K\in\mathcal{C}_1$ and 
$\theta_w \mathrm{A}I\in\mathcal{C}_2$. In particular, $\theta_w K$ is 
the maximal submodule of $\theta_w  I$, which belongs to $\mathcal{C}_1$.
Furthermore, the morphism  $\theta_w(\mathrm{nat}_I)$ is surjective.

At the same time, the module $\theta_w I$ is injective as $\theta_w$ is
right adjoint to the exact functor $\theta_{w^{-1}}$. From the definition
of $\mathrm{A}$ we have that the morphism $\mathrm{nat}_{\theta_w I}$
is surjective and that its kernel coincides with the maximal submodule 
of $\theta_w  I$, which belongs to $\mathcal{C}_1$. In other words, the
kernels of $\mathrm{nat}_{\theta_w I}$ and $\theta_w(\mathrm{nat}_{I})$
coincide.

Now the statement of the lemma follows from \cite[Lemma~1]{KM},
applied to the situation $\mathrm{F}=\mathrm{A}\theta_w$,
$\mathrm{G}=\theta_w\mathrm{A}$ and $\mathrm{H}=\theta_w$.
\end{proof}

Set $\overline{D}^{\hat{\mathbf{R}}}=\mathrm{A}L(\mathbf{w})$. 

\begin{lemma}\label{lem7}
\begin{enumerate}[(i)]
\item \label{lem7.1} $\overline{D}^{\hat{\mathbf{R}}}$ is
isomorphic to the maximal submodule of the module 
$P^{\hat{\mathbf{R}}}(\mathbf{w})$,
which contains the socle of $P^{\hat{\mathbf{R}}}(\mathbf{w})$ and such that
all other composition subquotients of $\overline{D}^{\hat{\mathbf{R}}}$ have
the  form $L(x)$, where $x<_{\mathtt{R}}\mathbf{w}$.
\item \label{lem7.2} We have $D^{\hat{\mathbf{R}}}\subset 
\overline{D}^{\hat{\mathbf{R}}}$ and the condition of Theorem~\ref{thm2}\eqref{thm2.2} is  equivalent to the equality
$D^{\hat{\mathbf{R}}}=\overline{D}^{\hat{\mathbf{R}}}$.
\end{enumerate}
\end{lemma}

\begin{proof}
As $P^{\hat{\mathbf{R}}}(\mathbf{w})$ is the injective envelope of
$L(\mathbf{w})$, the statement \eqref{lem7.1} follows immediately from
the definition of $\mathrm{A}$. 

The inclusion $D^{\hat{\mathbf{R}}}\subset \overline{D}^{\hat{\mathbf{R}}}$
follows from \cite[Lemmata~6-8]{MS2}. The rest of the statement \eqref{lem7.2} 
now follows from \eqref{lem7.1} and the definition of $\mathrm{A}$.
\end{proof}

\begin{lemma}\label{lemnew5}
For any $w\in W$ we have
\begin{displaymath}
\dim\mathrm{Hom}_{\mathfrak{g}}(L(\mathbf{w}),\theta_w L(\mathbf{w}))= 
\dim\mathrm{Hom}_{\mathfrak{g}}(\overline{D}^{\hat{\mathbf{R}}},
\theta_w \overline{D}^{\hat{\mathbf{R}}})
\end{displaymath}
\end{lemma}

\begin{proof}
Since $L(\mathbf{w})\in\mathcal{C}_2$, we have 
$\theta_w L(\mathbf{w})\in\mathcal{C}_2$ by Lemma~\ref{lemintro}. 
Hence, by the definition of $\mathrm{A}$, we have that $\mathrm{A}$ 
annihilates neither $L(\mathbf{w})$ nor any simple submodule of 
$\theta_w L(\mathbf{w})$. Applying $\mathrm{A}$ and using its 
definition we thus obtain an inclusion
\begin{displaymath}
\mathrm{Hom}_{\mathfrak{g}}(L(\mathbf{w}),\theta_w L(\mathbf{w}))\subset
\mathrm{Hom}_{\mathfrak{g}}(\mathrm{A}L(\mathbf{w}),
\mathrm{A}\theta_w L(\mathbf{w})).
\end{displaymath}
Using Lemma~\ref{lem4} and the definition of $\overline{D}^{\hat{\mathbf{R}}}$
we thus get the inclusion
\begin{equation}\label{eqnew5-101}
\mathrm{Hom}_{\mathfrak{g}}(L(\mathbf{w}),\theta_w L(\mathbf{w}))
\subset \mathrm{Hom}_{\mathfrak{g}}(\overline{D}^{\hat{\mathbf{R}}},
\theta_w \overline{D}^{\hat{\mathbf{R}}}).
\end{equation}

On the other hand, consider the short exact sequence
\begin{equation}\label{eqnew5-1}
0\to L(\mathbf{w})\to \overline{D}^{\hat{\mathbf{R}}}\to C\to 0,
\end{equation}
where $C$ is the cokernel. Applying the exact functor $\theta_w$ yields
the short exact sequence
\begin{equation}\label{eqnew5-2}
0\to \theta_w L(\mathbf{w})\to \theta_w \overline{D}^{\hat{\mathbf{R}}}
\to \theta_w C\to 0.
\end{equation}
Applying the bifunctor $\mathrm{Hom}_{\mathfrak{g}}({}_-,{}_-)$ from the
sequence \eqref{eqnew5-1} to the sequence \eqref{eqnew5-2} yields the
following commutative diagram with exact rows and columns
\begin{displaymath}
\xymatrix{
\mathrm{Hom}_{\mathfrak{g}}(C,\theta_w L(\mathbf{w}))
\ar@{^{(}->}[r]\ar@{^{(}->}[d] &
\mathrm{Hom}_{\mathfrak{g}}(C,\theta_w \overline{D}^{\hat{\mathbf{R}}})
\ar@{^{(}->}[d]\ar[r]&
\mathrm{Hom}_{\mathfrak{g}}(C,\theta_w C)\ar@{^{(}->}[d]\\
\mathrm{Hom}_{\mathfrak{g}}(\overline{D}^{\hat{\mathbf{R}}},
\theta_w L(\mathbf{w}))\ar@{^{(}->}[r]\ar[d] &
\mathrm{Hom}_{\mathfrak{g}}(\overline{D}^{\hat{\mathbf{R}}},
\theta_w \overline{D}^{\hat{\mathbf{R}}})\ar[r]\ar[d]&
\mathrm{Hom}_{\mathfrak{g}}(\overline{D}^{\hat{\mathbf{R}}},\theta_w C)\ar[d]\\
\mathrm{Hom}_{\mathfrak{g}}(L(\mathbf{w}),\theta_w L(\mathbf{w}))
\ar@{^{(}->}[r] &
\mathrm{Hom}_{\mathfrak{g}}(L(\mathbf{w}),
\theta_w\overline{D}^{\hat{\mathbf{R}}}) \ar[r]&
\mathrm{Hom}_{\mathfrak{g}}(L(\mathbf{w}),\theta_w C).\\
}
\end{displaymath}

We have  $C,\theta_w C\in\mathcal{C}_1$ by definitions and 
Lemma~\ref{lemintro}, and $L(\mathbf{w})\in\mathcal{C}_3$. This yields
$\mathrm{Hom}_{\mathfrak{g}}(L(\mathbf{w}),\theta_w C)=0$, which  implies  
\begin{displaymath}
\mathrm{Hom}_{\mathfrak{g}}(L(\mathbf{w}),\theta_w L(\mathbf{w}))=
\mathrm{Hom}_{\mathfrak{g}}(L(\mathbf{w}),
\theta_w\overline{D}^{\hat{\mathbf{R}}}).
\end{displaymath}

Since $C\in\mathcal{C}_1$ while $\overline{D}^{\hat{\mathbf{R}}},
\theta_w\overline{D}^{\hat{\mathbf{R}}}\in \mathcal{C}_2$ by definitions 
and  Lemma~\ref{lemintro}, we have 
$\mathrm{Hom}_{\mathfrak{g}}(C,\theta_w\overline{D}^{\hat{\mathbf{R}}})=0$,
which yields the inclusion
\begin{displaymath}
\mathrm{Hom}_{\mathfrak{g}}(\overline{D}^{\hat{\mathbf{R}}},
\theta_w \overline{D}^{\hat{\mathbf{R}}})\subset 
\mathrm{Hom}_{\mathfrak{g}}(L(\mathbf{w}),
\theta_w\overline{D}^{\hat{\mathbf{R}}}).
\end{displaymath}
The latter, together with the equality, obtained in the previous paragraph, 
implies the opposite to \eqref{eqnew5-101} inclusion
\begin{displaymath}
\mathrm{Hom}_{\mathfrak{g}}(\overline{D}^{\hat{\mathbf{R}}},
\theta_w \overline{D}^{\hat{\mathbf{R}}})\subset
\mathrm{Hom}_{\mathfrak{g}}(L(\mathbf{w}),\theta_w L(\mathbf{w}))
\end{displaymath}
and the statement of the lemma follows.
\end{proof}

\begin{lemma}\label{lemnew6}
The inclusion $L(\mathbf{w})\subset \overline{D}^{\hat{\mathbf{R}}}$
induces an isomorphism of $\mathfrak{g}$-bimodules as follows:
$\mathcal{L}(L(\mathbf{w}),L(\mathbf{w}))\cong
\mathcal{L}(\overline{D}^{\hat{\mathbf{R}}},
\overline{D}^{\hat{\mathbf{R}}})$.
\end{lemma}

\begin{proof}
Applying the bifunctor $\mathcal{L}({}_-,{}_-)$ to \eqref{eqnew5-1}
we get the following commutative diagram with exact rows and columns: 
\begin{equation}\label{eqnew6-1}
\xymatrix{
\mathcal{L}(C,L(\mathbf{w}))\ar@{^{(}->}[rr]\ar@{^{(}->}[d] &&
\mathcal{L}(C,\overline{D}^{\hat{\mathbf{R}}})\ar@{^{(}->}[d]\ar[rr]&&
\mathcal{L}(C,C)\ar@{^{(}->}[d]\\
\mathcal{L}(\overline{D}^{\hat{\mathbf{R}}},L(\mathbf{w}))
\ar@{^{(}->}[rr]\ar[d] &&
\mathcal{L}(\overline{D}^{\hat{\mathbf{R}}},
\overline{D}^{\hat{\mathbf{R}}})\ar[rr]\ar[d]&&
\mathcal{L}(\overline{D}^{\hat{\mathbf{R}}},C)\ar[d]\\
\mathcal{L}(L(\mathbf{w}),L(\mathbf{w}))\ar@{^{(}->}[rr] &&
\mathcal{L}(L(\mathbf{w}),\overline{D}^{\hat{\mathbf{R}}})\ar[rr]&&
\mathcal{L}(L(\mathbf{w}),C).\\
}
\end{equation}

Since for any $w\in W$ we have $C,\theta_w C\in\mathcal{C}_1$ by
definitions and Lemma~\ref{lemintro}, while 
$L(\mathbf{w})\in\mathcal{C}_3$, from \cite[6.8(3)]{Ja} we have
$\mathcal{L}(L(\mathbf{w}),C)=0$ implying
$\mathcal{L}(L(\mathbf{w}),L(\mathbf{w}))\cong
\mathcal{L}(L(\mathbf{w}),\overline{D}^{\hat{\mathbf{R}}})$.

Since for any $w\in W$ we have $\overline{D}^{\hat{\mathbf{R}}},
\theta_w \overline{D}^{\hat{\mathbf{R}}}\in\mathcal{C}_2$ by
definitions and Lemma~\ref{lemintro}, while 
$C\in\mathcal{C}_1$, from \cite[6.8(3)]{Ja} it follows
that $\mathcal{L}(C,\overline{D}^{\hat{\mathbf{R}}})=0$ implying 
$\mathcal{L}(\overline{D}^{\hat{\mathbf{R}}},
\overline{D}^{\hat{\mathbf{R}}})\subset
\mathcal{L}(L(\mathbf{w}),\overline{D}^{\hat{\mathbf{R}}})$.
Taking the above, Lemma~\ref{lemnew5} and \cite[6.8(3)]{Ja} into 
account yields $\mathcal{L}(L(\mathbf{w}),L(\mathbf{w}))\cong
\mathcal{L}(\overline{D}^{\hat{\mathbf{R}}},
\overline{D}^{\hat{\mathbf{R}}})$, which completes  the proof.
\end{proof}

\begin{proof}[Proof of the implication 
\eqref{thm2.2}$\Rightarrow$\eqref{thm2.1} in Theorem~\ref{thm2}.]
Because of the as\-sum\-pti\-on Theorem~\ref{thm2}\eqref{thm2.2},
from Lemma~\ref{lem7}\eqref{lem7.2} we have 
$\overline{D}^{\hat{\mathbf{R}}}={D}^{\hat{\mathbf{R}}}$.
The module $D^{\hat{\mathbf{R}}}$ is a quotients of the dominant Verma
module $\Delta(e)$ and hence $U(\mathfrak{g})$ surjects onto
$\mathcal{L}(D^{\hat{\mathbf{R}}},D^{\hat{\mathbf{R}}})$ by 
\cite[6.9(10)]{Ja}. Lemma~\ref{lemnew6} and the diagram  \eqref{eqnew6-1} 
now give the induced  surjection of $U(\mathfrak{g})$ onto 
$\mathcal{L}(L(\mathbf{w}),L(\mathbf{w}))$. This completes the proof.
\end{proof}

\begin{corollary}\label{cor14}
If the condition Theorem~\ref{thm2}\eqref{thm2.2} is satisfied, we
have the equality $\mathrm{Ann}_{U(\mathfrak{g})}(L(\mathbf{w}))=
\mathrm{Ann}_{U(\mathfrak{g})}(D^{\hat{\mathbf{R}}})$.
\end{corollary}

\begin{proof}
From $L(\mathbf{w})\subset D^{\hat{\mathbf{R}}}$ we have 
$\mathrm{Ann}_{U(\mathfrak{g})}(D^{\hat{\mathbf{R}}})\subset 
\mathrm{Ann}_{U(\mathfrak{g})}(L(\mathbf{w}))$. On the other hand,
from the previous proof we have 
\begin{multline*}
U(\mathfrak{g})/\mathrm{Ann}_{U(\mathfrak{g})}(D^{\hat{\mathbf{R}}})\cong
\mathcal{L}(D^{\hat{\mathbf{R}}},D^{\hat{\mathbf{R}}})\cong \\ \cong
\mathcal{L}(L(\mathbf{w}),L(\mathbf{w}))\cong
U(\mathfrak{g})/\mathrm{Ann}_{U(\mathfrak{g})}(L(\mathbf{w})),
\end{multline*}
which implies the statement.
\end{proof}

\begin{lemma}\label{lem25}
\begin{enumerate}[(i)]
\item\label{lem25.1} $\mathrm{A}D^{\hat{\mathbf{R}}}\cong
\overline{D}^{\hat{\mathbf{R}}}$.
\item\label{lem25.2} $\mathrm{A}\overline{D}^{\hat{\mathbf{R}}}\cong
\overline{D}^{\hat{\mathbf{R}}}$.
\item\label{lem25.3} For any $w\in W$ there is an isomorphism
\begin{displaymath}
\mathrm{Hom}_{\mathfrak{g}}(D^{\hat{\mathbf{R}}},\theta_w
\overline{D}^{\hat{\mathbf{R}}}) \cong
\mathrm{Hom}_{\mathfrak{g}}(\overline{D}^{\hat{\mathbf{R}}},\theta_w
\overline{D}^{\hat{\mathbf{R}}}).
\end{displaymath}
\end{enumerate}
\end{lemma}

\begin{proof}
Consider the short exact sequence
\begin{equation}\label{eq26}
0\to  D^{\hat{\mathbf{R}}}\to \overline{D}^{\hat{\mathbf{R}}}\to
C\to 0,
\end{equation} 
where $C\in\mathcal{C}_1$ is the cokernel. From the definition of 
$\mathrm{A}$ we have $\mathrm{A}C=0$. Applying now $\mathrm{A}$
to \eqref{eq26} and using the left exactness of $\mathrm{A}$ yields
the statement \eqref{lem25.1}. The statement \eqref{lem25.2} follows
immediately from the definition of $\mathrm{A}$. 

Since $C\in \mathcal{C}_1$ and $\overline{D}^{\hat{\mathbf{R}}},
\theta_w\overline{D}^{\hat{\mathbf{R}}}\in \mathcal{C}_2$, applying
$\mathrm{Hom}_{\mathfrak{g}}({}_-,\theta_w
\overline{D}^{\hat{\mathbf{R}}})$ to \eqref{eq26} yields the
inclusion
\begin{equation}\label{eq27}
\mathrm{Hom}_{\mathfrak{g}}(\overline{D}^{\hat{\mathbf{R}}},\theta_w
\overline{D}^{\hat{\mathbf{R}}})\subset
\mathrm{Hom}_{\mathfrak{g}}(D^{\hat{\mathbf{R}}},\theta_w
\overline{D}^{\hat{\mathbf{R}}}).
\end{equation}
On the other hand, the functor $\mathrm{A}$ annihilates neither the
socle of ${D}^{\hat{\mathbf{R}}}$ nor any submodule in the socle of 
$\theta_w \overline{D}^{\hat{\mathbf{R}}}$. Hence from the definition
of $\mathrm{A}$ we have the inclusion
\begin{displaymath}
\mathrm{Hom}_{\mathfrak{g}}(D^{\hat{\mathbf{R}}},\theta_w
\overline{D}^{\hat{\mathbf{R}}})\subset
\mathrm{Hom}_{\mathfrak{g}}(\mathrm{A} D^{\hat{\mathbf{R}}},
\mathrm{A} \theta_w \overline{D}^{\hat{\mathbf{R}}}).
\end{displaymath}
Using \eqref{lem25.1}, Lemma~\ref{lem4} and \eqref{lem25.2} 
we obtain 
\begin{displaymath}
\mathrm{Hom}_{\mathfrak{g}}(\mathrm{A} D^{\hat{\mathbf{R}}},
\mathrm{A}\theta_w \overline{D}^{\hat{\mathbf{R}}})=
\mathrm{Hom}_{\mathfrak{g}}(\overline{D}^{\hat{\mathbf{R}}},
\theta_w\mathrm{A} \overline{D}^{\hat{\mathbf{R}}})=
\mathrm{Hom}_{\mathfrak{g}}(\overline{D}^{\hat{\mathbf{R}}},\theta_w
\overline{D}^{\hat{\mathbf{R}}}),
\end{displaymath}
which implies that the inclusion \eqref{eq27} is in fact an isomorphism.
This completes the proof.
\end{proof}

\begin{proof}[Proof of the implication 
\eqref{thm2.1}$\Rightarrow$\eqref{thm2.2} in Theorem~\ref{thm2}.]
The inclusion $L(\mathbf{w})\subset D^{\hat{\mathbf{R}}}$ induces
the inclusion $\mathrm{Ann}_{U(\mathfrak{g})}(D^{\hat{\mathbf{R}}})\subset
\mathrm{Ann}_{U(\mathfrak{g})}(L(\mathbf{w}))$, which, in turn, induces
the surjection
\begin{equation}\label{eq20}
U(\mathfrak{g})/\mathrm{Ann}_{U(\mathfrak{g})}(D^{\hat{\mathbf{R}}})
\tto  
U(\mathfrak{g})/\mathrm{Ann}_{U(\mathfrak{g})}(L(\mathbf{w})).
\end{equation}

Assume that the condition of Theorem~\ref{thm2}\eqref{thm2.2} is not 
satisfied. As we have 
$\mathcal{L}(D^{\hat{\mathbf{R}}},D^{\hat{\mathbf{R}}})\cong
U(\mathfrak{g})/\mathrm{Ann}_{U(\mathfrak{g})}(D^{\hat{\mathbf{R}}})$
by \cite[6.9(10)]{Ja}, from the latter formula and \eqref{eq20} it follows
that the inequality
\begin{equation}\label{eq21}
\mathcal{L}(D^{\hat{\mathbf{R}}},D^{\hat{\mathbf{R}}})\subsetneq 
\mathcal{L}(L(\mathbf{w}),L(\mathbf{w}))\cong
\mathcal{L}(\overline{D}^{\hat{\mathbf{R}}},\overline{D}^{\hat{\mathbf{R}}}).
\end{equation}
would imply that the algebra $U(\mathfrak{g})$ does not surjects onto
$\mathcal{L}(L(\mathbf{w}),L(\mathbf{w}))$. Hence we are now left to prove
the inequality \eqref{eq21}.

We apply the bifunctor $\mathcal{L}({}_-,{}_-)$ to the short exact
sequence \eqref{eq26}, where the cokernel $C\neq 0$ by 
Lemma~\ref{lem7}\eqref{lem7.2}. Since $C\in\mathcal{C}_1$ and 
$D^{\hat{\mathbf{R}}}$,
$\theta_w D^{\hat{\mathbf{R}}}$, $\overline{D}^{\hat{\mathbf{R}}}$ and
$\theta_w \overline{D}^{\hat{\mathbf{R}}}$ are in $\mathcal{C}_2$
for all $w\in W$, by \cite[6.8(3)]{Ja} we obtain
the following commutative diagram with exact rows and columns: 
\begin{equation}\label{eq23}
\xymatrix{
0\ar[rr]\ar[d] && 0\ar[d]\ar[rr]&&
\mathcal{L}(C,C)\ar@{^{(}->}[d]\\
\mathcal{L}(\overline{D}^{\hat{\mathbf{R}}},D^{\hat{\mathbf{R}}})
\ar@{^{(}->}[rr]\ar@{^{(}->}[d] &&
\mathcal{L}(\overline{D}^{\hat{\mathbf{R}}},
\overline{D}^{\hat{\mathbf{R}}})\ar[rr]\ar[d]^{\wr}&&
\mathcal{L}(\overline{D}^{\hat{\mathbf{R}}},C)\ar[d]\\
\mathcal{L}(D^{\hat{\mathbf{R}}},D^{\hat{\mathbf{R}}})\ar@{^{(}->}[rr] &&
\mathcal{L}(D^{\hat{\mathbf{R}}},\overline{D}^{\hat{\mathbf{R}}})
\ar[rr]^{\alpha}&& \mathcal{L}(D^{\hat{\mathbf{R}}},C),\\
}
\end{equation}
where the isomorphism in the second column follows from 
Lemma~\ref{lem25}\eqref{lem25.3}. To complete the proof it is thus enough to
show that the map $\alpha$ on the diagram \eqref{eq23} is non-zero.

Pick some simple submodule $L(x)\subset C$ (recall once more that
$C\neq 0$ by Lemma~\ref{lem7}\eqref{lem7.2}). Then, using the adjointness
and defining properties of  projective functors, we have
\begin{equation}\label{eq29}
\begin{array}{rcl}
\mathbb{C}&=&\mathrm{Hom}_{\mathfrak{g}}(P^{\hat{\mathbf{R}}}(x),L(x)) \\
&=&\mathrm{Hom}_{\mathfrak{g}}(\theta_x P^{\hat{\mathbf{R}}}(e),L(x)) \\
&=&\mathrm{Hom}_{\mathfrak{g}}(P^{\hat{\mathbf{R}}}(e),\theta_{x^{-1}}L(x))\\
&\subset&\mathrm{Hom}_{\mathfrak{g}}(P^{\hat{\mathbf{R}}}(e),
\theta_{x^{-1}}C).
\end{array}
\end{equation}

Applying the bifunctor $\mathrm{Hom}_{\mathfrak{g}}({}_-,{}_-)$
from the short exact sequence
\begin{displaymath}
0\to K\to  P^{\hat{\mathbf{R}}}(e)\to D^{\hat{\mathbf{R}}}\to 0,
\end{displaymath}
where $K$ is just the kernel of the natural projection
$P^{\hat{\mathbf{R}}}(e)\tto D^{\hat{\mathbf{R}}}$
(note that $K\in\mathcal{C}_1$ by \cite[Lemmata 6-8]{MS2}), to 
the short exact sequence
\begin{displaymath}
0\to  \theta_{x^{-1}} D^{\hat{\mathbf{R}}}\to
\theta_{x^{-1}}\overline{D}^{\hat{\mathbf{R}}}\to
\theta_{x^{-1}} C\to 0
\end{displaymath}
we obtain the following commutative diagram with exact rows and columns:
{\small
\begin{displaymath}
\xymatrix{ 
\mathrm{Hom}_{\mathfrak{g}}(D^{\hat{\mathbf{R}}},\theta_{x^{-1}}
D^{\hat{\mathbf{R}}})\ar@{^{(}->}[r]\ar[d]^{\wr}
&\mathrm{Hom}_{\mathfrak{g}}(D^{\hat{\mathbf{R}}},\theta_{x^{-1}}
\overline{D}^{\hat{\mathbf{R}}})\ar[r]\ar[d]^{\wr}\ar@{:>}[dr]^{\beta}
&\mathrm{Hom}_{\mathfrak{g}}(D^{\hat{\mathbf{R}}},\theta_{x^{-1}}
C)\ar@{^{(}->}[d]\\
\mathrm{Hom}_{\mathfrak{g}}(P^{\hat{\mathbf{R}}}(e),\theta_{x^{-1}}
D^{\hat{\mathbf{R}}})\ar@{^{(}->}[r]\ar[d]
&\mathrm{Hom}_{\mathfrak{g}}(P^{\hat{\mathbf{R}}}(e),\theta_{x^{-1}}
\overline{D}^{\hat{\mathbf{R}}})\ar@{->>}[r]\ar[d]
&\mathrm{Hom}_{\mathfrak{g}}(P^{\hat{\mathbf{R}}}(e),\theta_{x^{-1}}
C)\ar[d]\\
0\ar[r]
&0\ar[r]
&\mathrm{Hom}_{\mathfrak{g}}(K,\theta_{x^{-1}}
C),\\
}
\end{displaymath}
}
where the second row is exact as $P^{\hat{\mathbf{R}}}(e)$ is
projective in $\mathcal{O}_0^{\hat{\mathbf{R}}}$, and the zeros in the 
third row follow from the fact that $K\in\mathcal{C}_1$ 
while $\theta_{x^{-1}}\overline{D}^{\hat{\mathbf{R}}} \in\mathcal{C}_2$.
From \eqref{eq29} it follows that the composition $\beta$ is a surjection
onto a non-zero vector space, hence is a non-zero map. From the definitions 
we have that $\beta\neq 0$ implies $\alpha\neq 0$. This completes the proof.
\end{proof}

\subsection{A sufficient condition for the negative answer}\label{s2.2}

Let $\Lambda$ be the basic finite-dimensional associative algebra,
whose module category is equivalent to $\mathcal{O}_0$. The algebra
$\Lambda$ is Koszul (see \cite{So}) so we can fix the positive Koszul
$\mathbb{Z}$-grading on $\Lambda$. Let $\Lambda\mathrm{-gmod}$ denote 
the category of finite-dimensional graded $\Lambda$-modules. For $x\in 
\hat{\mathbf{R}}$ let $\mathtt{P}^{\hat{\mathbf{R}}}(x)$ denote the standard 
graded lift of ${P}^{\hat{\mathbf{R}}}(x)$ with head concentrated in 
degree zero (see \cite[4.3]{MS1}), and $\mathtt{L}(x)$ denote the standard
graded lift of the corresponding simple quotient (concentrated in degree 
zero). For $w\in W$ we denote by $\hat{\theta}_w$ the standard graded lifts 
of the functors $\theta_w$, see \cite[Section~8]{St}. Finally, let 
$\mathbf{a}:W\to\mathbb{Z}$ denote Lusztig's $\mathbf{a}$-function
(see \cite{Lu}), which is uniquely determined by the properties that it 
is constant on the two-sided cells of $W$ and equals the length of
$w'_0$ on every $w'_0$, which is the longest element of a parabolic 
subgroup  of $W$. 

If $\mathtt{M}$ is a graded module, then $\mathtt{M}=
\oplus_{i\in\mathbb{Z}} \mathtt{M}_i$ is the decomposition of $\mathtt{M}$
into a direct sum of graded components. As usually, for $k\in\mathbb{Z}$
we denote by $\langle k\rangle:\Lambda\mathrm{-gmod}\to\Lambda\mathrm{-gmod}$ 
the functor, which shifts the grading such that 
$\mathtt{M}\langle k\rangle_{i}=\mathtt{M}_{i+k}$.

\begin{lemma}\label{lem554}
Let $\mathbf{w}\in W$ be an involution and $\mathtt{M}=
\hat{\theta}_{\mathbf{w}}L(\mathbf{w})$. Then:
\begin{enumerate}[(i)]
\item\label{lem554.1} $\mathtt{M}_i=0$ for all $i$ such that 
$|i|>\mathbf{a}(\mathbf{w})$. 
\item\label{lem554.2} $\mathtt{M}_{\mathbf{a}(\mathbf{w})}$ is
the simple socle of $\mathtt{M}$ (which is isomorphic to
the module $\mathtt{L}(\mathbf{w})\langle -\mathbf{a}(\mathbf{w})\rangle$).
\end{enumerate}
\end{lemma}

\begin{proof}
Since $\mathbf{a}$ is an invariant of two-sided cells, by
\cite[Theorem~18]{MS1} we may without loss of generality assume that
$\mathbf{w}$ is the maximal element of some parabolic subgroup. 
For such $\mathbf{w}$ the statement \eqref{lem554.1} follows 
immediately from \cite[Theorem~8.2]{St}. Moreover, the same argument
implies $\mathtt{M}_{\mathbf{a}(\mathbf{w})}\neq 0$.

As $\Lambda$ is positively graded and $\mathtt{M}$ is injective
(the latter follows from \cite[Section~5]{MS1} and 
\cite[Key statement]{MS2}), $\mathtt{M}_{\mathbf{a}(\mathbf{w})}\neq 0$
must be the simple socle of $\mathtt{M}$. This completes the proof.
\end{proof}

\begin{theorem}\label{thm555}
Let $\mathbf{w}\in W$ be an involution and $\mathtt{M}=
\hat{\theta}_{\mathbf{w}}\mathtt{L}(\mathbf{w})$. Assume that there exists
$x\in W$ such that $x<_{\mathtt{R}}\mathbf{w}$ and 
\begin{displaymath}
[\mathtt{M}:\mathtt{L}(x)\langle 1-\mathbf{a}(\mathbf{w})\rangle]>
[\mathtt{P}^{\hat{\mathbf{R}}}(e):\mathtt{L}(x)\langle 1-\mathbf{a}(\mathbf{w})\rangle]. 
\end{displaymath}
Then Kostant's problem has the negative answer for $L(\mathbf{w})$.
\end{theorem}

\begin{proof}
Let $\mathtt{N}$ be the quotient of $\mathtt{M}$ modulo $D^{\hat{\mathbf{R}}}$.
As $D^{\hat{\mathbf{R}}}$ is non-zero, it must contain the socle of
$\mathtt{M}$. Hence $\mathtt{N}_i=0$ for all $i\geq \mathbf{a}(\mathbf{w})$
by Lemma~\ref{lem554}. By our assumption, 
$\mathtt{N}_{\mathbf{a}(\mathbf{w})-1}$ contains at least one
copy of $\mathtt{L}(x)\langle 1-\mathbf{a}(\mathbf{w})\rangle$.

Since $\Lambda$ is positively graded and $\mathtt{N}_i=0$ for all 
$i\geq \mathbf{a}(\mathbf{w})$, the space
$\mathtt{N}_{\mathbf{a}(\mathbf{w})-1}$ belongs to the socle of $\mathtt{N}$.
Thus the condition of Theorem~\ref{thm2}\eqref{thm2.2} is not satisfied
and the answer to Kostant's problem for $L(\mathbf{w})$ is negative 
by Theorem~\ref{thm2}.
\end{proof}

\begin{remark}\label{rem556}
{\rm  
As $\mathtt{P}^{\hat{\mathbf{R}}}(e)$ is a quotient of the graded
dominant Verma module $\Delta(e)$, in Lemma~\ref{thm555} one could
use a stronger assumption
\begin{displaymath}
[\mathtt{M}:\mathtt{L}(x)\langle 1-\mathbf{a}(\mathbf{w})\rangle]>
[\Delta(e):\mathtt{L}(x)\langle 1-\mathbf{a}(\mathbf{w})\rangle] 
\end{displaymath}
with the same result.
}
\end{remark}

\begin{remark}\label{rem557}
{\rm  
The numerical condition of Theorem~\ref{thm555} is relatively easy to
check for example using the computer, because it can be easily formulated
in terms of Kazhdan-Lusztig combinatorics, \cite{KL,BB}. 
Via the standard categorification approach to $\mathcal{O}$
(see for example \cite[3.4]{MS1}), the characters of graded 
$\Lambda$-modules can be considered as elements of the Hecke algebra $\mathcal{H}$ of $W$ (such that Verma modules correspond to the standrad
basis of $\mathcal{H}$, projective modules correspond to the Kazhdan-Lusztig
basis, and simple modules correspond to the dual Kazhdan-Lusztig basis). 
There are effective algorithms, which allow one to multiply elements
of $\mathcal{H}$ and to transform them from one of the mentioned basis
to the other. Some of the applications, presented in the next section 
are obtained using this approach.
}
\end{remark}

\begin{remark}\label{rem558}
{\rm 
The statement of Lemma~\ref{lem554} has a strong resemblance with
\cite[Theorem~16]{Ma2}, and is in some sense the Koszul dual of it
(see the proof of \cite[Theorem~16]{Ma2} for details).
}
\end{remark}

\section{Applications}\label{s3}

\subsection{Kostant's problem for the socle of the dominant Verma module
in a parabolic category}\label{s3.1}

Let $\mathfrak{p}\subset \mathfrak{g}$ be a parabolic subcategory
containing $\mathfrak{h}\oplus\mathfrak{n}_+$, and 
$\mathcal{O}_0^{\mathfrak{p}}$ be the corresponding parabolic subcategory
of $\mathcal{O}_0$ in the sense of \cite{RC}. Let $W'\subset W$ be the
Weyl group of the Levi factor of $\mathfrak{p}$, $w_0$ be the longest
element in $W$ and $w'_0$ be the longest element in $W'$. Then 
$\mathcal{O}_0^{\mathfrak{p}}=\mathcal{O}_0^{\hat{\mathbf{R}}}$,
where $\mathbf{R}$ is the right cell of the element $w'_0w_0$,
see \cite[Remark~14]{MS1}. Let $\mathbf{w}$ be the involution in 
$\mathbf{R}$.

\begin{corollary}\label{cor31}
Kostant's problem has the positive answer for $L(\mathbf{w})$. 
\end{corollary}

\begin{proof}
The category  $\mathcal{O}_0^{\mathfrak{p}}$ is known to be a highest
weight category in the sense of \cite{CPS}. Thus any projective-injective
module in $\mathcal{O}_0^{\mathfrak{p}}$ is tilting in the sense of
\cite{Ri}, in particular, it has a filtration by standard modules
(i.e. generalized Verma modules, induced from simple finite-dimensional
$\mathfrak{p}$-modules). In particular, the dominant standard module
$P^{\hat{\mathbf{R}}}(e)$ is a submodule of 
$P^{\hat{\mathbf{R}}}(\mathbf{w})$, and the cokernel
of this inclusion again has a filtration by standard modules.
Since all standard modules belong to $\mathcal{C}_2$ by \cite{Ir}
(see also \cite[Theorem~5.1]{MS3} for a short argument), we obtain that
the condition of Theorem~\ref{thm2}\eqref{thm2.2} is satisfied and
hence Kostant's problem has the positive answer for $L(\mathbf{w})$ 
by Theorem~\ref{thm2}.
\end{proof}

\begin{remark}\label{rem3101}
{\rm  
The fact $P^{\hat{\mathbf{R}}}(e)\subset
P^{\hat{\mathbf{R}}}(\mathbf{w})$ implies that the statement of 
Conjecture~\ref{conj1} is true if $\mathbf{R}$ contains some
$w'_0w_0$.
}
\end{remark}

\subsection{Kostant's problem for $L(s)$, where $s$ is a simple reflection}\label{s3.2}

\begin{corollary}[\cite{Ma}]\label{cor32}
Let $s\in W$ be a simple reflection. Then
Kos\-tant's problem has the positive answer for $L(s)$. 
\end{corollary}

\begin{proof}
The only element of $W$, which is strictly smaller than $s$ with respect
to the order $<_{\mathtt{R}}$ is the identity element $e$. As, by 
adjointness,
\begin{displaymath}
\dim\mathrm{Hom}_{\mathfrak{g}}(P^{\hat{\mathbf{R}}}(e),\theta_s L(s))= 
\dim\mathrm{Hom}_{\mathfrak{g}}(P^{\hat{\mathbf{R}}}(s),L(s))=1,
\end{displaymath}
the module $L(e)$ occurs in $\theta_s L(s)$ with multiplicity one, and hence
$L(e)$ does not occur in  the cokernel of the inclusion
$D^{\hat{\mathbf{R}}}\subset \theta_s L(s)$ at all. Therefore the condition 
of  Theorem~\ref{thm2}\eqref{thm2.2} is obviously satisfied and
hence Kostant's problem has the positive answer for $L(s)$ 
by Theorem~\ref{thm2}.
\end{proof}

\subsection{Kostant's problem for $L(st)$, where $s$ and $t$ are
commuting simple reflections}\label{s3.3}

Here we generalize the counterexample, constructed in \cite[Section~5]{MS2}.
Let $s_i=(i,i+1)$, $i=1,\dots,n-1$, be the $i$-th simple reflection in $W$.
We recall that for a simple reflection $s\in W$ and any $x\in W$ such that
$xs<x$ with respect to the Bruhat order we have that the module
$\hat{\theta}_s\mathtt{L}(x)$ is self-dual with simple head and socle
and we moreover have the following graded picture (the middle row of which 
is in degree $0$):
\begin{equation}\label{s3.3-eq1}
\xymatrix@!=0.6pc{ 
&&&\mathtt{L}(x)\langle 1\rangle\ar[rd]\ar[ld]&\\
\hat{\theta}_s\mathtt{L}(x):&&\mathtt{L}(xs)\ar[rd]&&X\ar[ld]\\
&&&\mathtt{L}(x)\langle -1\rangle,&
}
\end{equation}
where $X$ is a direct sum of $L(y)$'s such that $ys>y$ with multiplicity
$\mu(x,y)$, where $\mu$ is Kazhdan-Lusztig's $\mu$-function, \cite{KL}.
The formula \eqref{s3.3-eq1} is a standard corollary of (now proved) 
Kazhdan-Lusztig's conjecture in equivalent Vogan's form (see 
\cite{KL,GJ0,Vo}). We also refer to Remark~\ref{rem557}  and to \cite[Section~8]{St} for the appropriate graded reformulation.
The arrows on \eqref{s3.3-eq1} schematically represent the action of 
the algebra $\Lambda$.

\begin{corollary}\label{cor33}
Let $s_i$ and $s_j$ be two commuting different simple reflections in $W$
(i.e. $|i-j|>1$).  Then Kostant's problem has the positive answer for 
$L(s_is_j)$ if and only if $|i-j|>2$. 
\end{corollary}

\begin{proof}
Without loss of generality we assume  $j>i$. Let $\mathbf{R}_e=\{e\}$,
$\mathbf{R}_i$ denote the right cell of $s_i$, $\mathbf{R}_j$ denote the 
right cell of $s_j$, and $\mathbf{R}$ denote the right cell of $s_is_j$.
Then the Hasse diagram of $<_\mathtt{R}$ on the set 
$\{\mathbf{R}_e,\mathbf{R}_i,\mathbf{R}_j,\mathbf{R}\}$, where
$\mathbf{R}$ is the maximum element, is as follows:
\begin{displaymath}
\xymatrix@!=0.9pc{
&\mathbf{R}\ar@{-}[ld]\ar@{-}[rd]&\\ 
\mathbf{R}_i\ar@{-}[rd]&&\mathbf{R}_j\ar@{-}[ld]\\
&\mathbf{R}_e,&\\
}
\end{displaymath}
and we further have
\begin{gather*}
\mathbf{R}_i=\{s_i,s_is_{i-1},\dots,s_is_{i-1}\dots 
s_{1},s_is_{i+1},\dots, s_is_{i+1}\dots s_{n-1}\};\\
\mathbf{R}_j=\{s_j,s_js_{j-1},\dots,s_js_{j-1}\dots 
s_{1},s_js_{j+1},\dots, s_js_{j+1}\dots s_{n-1}\}.
\end{gather*}
A direct calculation gives $\theta_{s_i}\theta_{s_j}=\theta_{s_is_j}=
\theta_{s_j}\theta_{s_i}$.

Assume first that $j=i+2$. Since both $s_is_{i+2}$ and $s_is_{i+1}s_{i+2}$
are Boolean elements of $W$ (in the sense of \cite{Mm}), we have that
the Kazhdan-Lusztig polynomial $P_{s_is_{i+2},s_is_{i+1}s_{i+2}}(q)=1$ by
\cite[Theorem~5.4]{Mm} and hence $\mu(s_is_{i+2},s_is_{i+1}s_{i+2})=1$ as 
well by definition. This yields that $\mathrm{Ext}^1_{\mathcal{O}}
(L(s_is_{i+2}),L(s_is_{i+1}s_{i+2}))\neq 0$ and thus $L(s_is_{i+1}s_{i+2})$
occurs as a composition subquotient in $\theta_{s_i}L(s_is_{i+2})$
(as a direct summand of $X$ in \eqref{s3.3-eq1}). Applying 
\eqref{s3.3-eq1} we get that 
$\mathtt{L}(s_is_{i+1}s_{i+2})\langle -1\rangle$
occurs as a composition subquotient in 
$\hat{\theta}_{s_is_{i+2}}\mathtt{L}(s_is_{i+2})$. Note that we have
$s_is_{i+1}s_{i+2}<_{\mathtt{R}}s_is_{i+2}$. At the same time from 
\cite[Lemma~7.2.5]{Di} it follows that $\mathtt{P}^{\hat{\mathbf{R}}}(e)_1$
contains only composition subquotients of the form 
$\mathtt{L}(s_k)\langle -1\rangle$, $k=1,\dots,n-1$. Hence the 
numerical assumption of Theorem~\ref{thm555} is satisfied and 
therefore the answer to Kostant's problem for $L(s_is_{i+2})$ is 
negative by Theorem~\ref{thm555}.

If $j>i+2$, a similar application of \cite[Theorem~5.4]{Mm} yields
$\mu(s_is_j,s_is_{i+1}\dots s_{j-1}s_j)=0$ and also
$\mu(s_is_j,s_js_{j-1}\dots s_{i+1}s_i)=0$. The only other elements of
$\mathbf{R}_i$ and $\mathbf{R}_j$, comparable with $s_is_j$ with respect
to the Bruhat order, are $s_i$ and $s_j$ respectively. Because of
\eqref{s3.3-eq1}, this means that 
\begin{displaymath}
\xymatrix@!=0.6pc{ 
&&&\mathtt{L}(s_is_j)\langle 1\rangle\ar[rd]\ar[ld]&\\
\hat{\theta}_{s_i}\mathtt{L}(s_is_j):&&\mathtt{L}(s_j)\ar[rd]&&X\ar[ld]\\
&&&\mathtt{L}(s_is_j)\langle -1\rangle,&
}
\end{displaymath}
where $X$ is a direct sum of simple modules $L(y)$, $y\in\mathbf{R}$.
Applying now $\hat{\theta}_{s_j}$ and using \eqref{s3.3-eq1} again we
obtain the following graded filtration for the module 
$\hat{\theta}_{s_is_j}\mathtt{L}(s_is_j)$:
\begin{equation}\label{s3.3-eq2}
\xymatrix@!=0.6pc{ 
&&&&\mathtt{L}(s_is_j)\langle 2\rangle\ar[dllll]\ar[d]\ar[drr]\ar[drrrr]&&&&\\
\mathtt{L}(s_j)\langle 1\rangle\ar[d]\ar[drr]\ar[drrrrrr]&&&&
\mathtt{L}(s_i)\langle 1\rangle\ar[d]&&
Y\langle 1\rangle\ar[dll]&&X'\langle 1\rangle\ar[d]\ar[dll]\\
\mathtt{L}(e)\ar[d]&&Z\ar[dll]&&\mathtt{L}(s_is_j)\ar[dllll]\ar[drrrr]&&
\mathtt{L}(s_is_j)\ar[dll]\ar[d]&&U\ar[d]\\
\mathtt{L}(s_j)\langle -1\rangle\ar[drrrr]&&&&
\mathtt{L}(s_i)\langle -1\rangle\ar[d]&&
Y\langle -1\rangle\ar[dll]&&X'\langle -1\rangle\ar[dllll]\\
&&&&\mathtt{L}(s_is_j)\langle -2\rangle,&&&&\\
}
\end{equation}
where $Z$ is a direct sum of simples modules of the form $L(y)$, 
$y\in \mathbf{R}_j$; $Y$ is a direct sum of simples modules of the 
form $L(y)$, $y\in \mathbf{R}$; and $X'$ is a direct summand of $X$.
Note that the arrows on \eqref{s3.3-eq2} (which are supposed to 
schematically represent the action of $\Lambda$) show only the
part of the action, which obviously comes from \eqref{s3.3-eq1}, but 
they do not show the whole action. From \cite[Lemma~7.2.5]{Di} it
follows that the module $\mathtt{D}^{\hat{\mathbf{R}}}$ looks as
follows:
\begin{displaymath}
\xymatrix@!=0.6pc{ 
\mathtt{D}^{\hat{\mathbf{R}}}:&&&\mathtt{L}(e)\ar[dr]\ar[dl]&\\
&&\mathtt{L}(s_i)\langle -1\rangle\ar[dr]& 
&\mathtt{L}(s_j)\langle -1\rangle\ar[dl]\\
&&&\mathtt{L}(s_is_j)\langle -2\rangle&
}
\end{displaymath}
Now we have to analyze \eqref{s3.3-eq2} to determine the cokernel $C$
of the inclusion $\mathtt{D}^{\hat{\mathbf{R}}}\subset \hat{\theta}_{s_is_j}\mathtt{L}(s_is_j)$. $C$ obviously contains both
$Y\langle -1\rangle$ and $X'\langle -1\rangle$, but all direct summands
of these modules have the form $L(y)$, $y\in\mathbf{R}$, by above.
None of the simple subquotients of $U$ can be in $C$ by \eqref{s3.3-eq1}.
Similarly one excludes $\mathtt{L}(s_i)\langle 1\rangle$ and
$\mathtt{L}(s_i)\langle 1\rangle$. All simple submodules in $Z$ have the form
$L(y)$, $y\in\mathbf{R}_j$. Considering 
$\hat{\theta}_{s_is_j}\mathtt{L}(s_is_j)=
\hat{\theta}_{s_i}\hat{\theta}_{s_j}\mathtt{L}(s_is_j)$ and using
the same arguments as above one shows that none of the simple submodules
of $Z$ belongs to $C$. Hence $C$ contains only simple modules of the form 
$L(y)$, $y\in\mathbf{R}$. Thus the condition of 
Theorem~\ref{thm2}\eqref{thm2.2} is satisfied and therefore
Kostant's problem has the positive answer
for $L(s_is_j)$ by Theorem~\ref{thm2}. This completes the proof.
\end{proof}

\subsection{Kostant's problem for $\mathfrak{sl}_n$, $n\leq 3$}\label{s3.4}

\begin{proposition}\label{prop34}
Assume that $n\leq 3$ and $w\in W$. Then Kostant's problem has the 
positive answer for $L(w)$. 
\end{proposition}

\begin{proof}
The statement is trivial for $n=1$. In the case $n=2$ for $w=e$ the
statement follows from \cite[6.9(10)]{Ja} 
(as $L(e)$ is a quotient of the dominant Verma module)
and for $w=s_1$ it follows from \cite[Corollary~6.4]{Jo2}
(as $L(s_1)$ is a Verma module).

Finally, in the case $n=3$ for $w=e$ the statement follows, as above, from
\cite[6.9(10)]{Ja}, for $w=s_1,s_2$ it follows from 
Corollary~\ref{cor32}, for $w=s_1s_1,s_2s_1$ it follows from 
\cite[Theorem~4.4]{GJ}, and, finally, for $w=s_1s_2s_1$ it follows, 
as above, from \cite[Corollary~6.4]{Jo2}.
\end{proof}

\subsection{Kostant's problem for  $\mathfrak{sl}_4$}\label{s3.5}

\begin{proposition}\label{prop35}
Assume that $n=4$ and $w\in W$. Then Kos\-tant's problem has the 
positive answer for $L(w)$ if and only if $w\neq s_1s_3,s_2s_1s_3$. 
\end{proposition}

\begin{proof}
The group $S_4$ has 10 involutions: $e$, $s_1$, $s_2$, $s_3$,
$s_1s_3$, $s_1s_2s_1$, $s_3s_2s_3$, $s_2s_1s_3s_2$,
$s_1s_2s_3s_2s_1$, and $s_2s_1s_2s_3s_2s_1$. The module $L(e)$ 
is a quotient of the dominant Verma module, and hence for $L(e)$
the claim follows from \cite[6.9(10)]{Ja}. The module 
$L(s_2s_1s_2s_3s_2s_1)$ is a Verma  module and hence for this module
the claim follows  from \cite[Corollary~6.4]{Jo2}. For $L(s_1)$, 
$L(s_2)$, $L(s_3)$ the claim follows from Corollary~\ref{cor32}.
The left cell of each of the elements $s_1s_2s_1$, $s_3s_2s_3$, 
$s_2s_1s_3s_2$, $s_1s_2s_3s_2s_1$ contains an element of the form
$w'_0w_0$, where $w'_0$ is the longest element of some parabolic 
subgroup. Hence for $L(s_1s_2s_1)$, $L(s_3s_2s_3)$, $L(s_2s_1s_3s_2)$ and
$L(s_1s_2s_3s_2s_1)$ the claim follows from \cite[Theorem~4.4]{GJ}
and \cite[Theorem~60]{MS1}.  Finally, for $L(s_1s_3)$ the claim follows
from Corollary~\ref{cor33} (or \cite[Theorem~13]{MS2}). Note that the
answer is negative only in the case of $L(s_1s_3)$. The left cell
of $s_1s_3$ contains one more element, namely $s_2s_1s_3$. The statement
of the proposition now follows from \cite[Theorem~60]{MS1}.
\end{proof}

\subsection{Kostant's problem for  $\mathfrak{sl}_5$}\label{s4.6}

\begin{proposition}\label{prop22}
Assume that $n=5$ and $w\in W$. Then Kos\-tant's problem has the 
positive answer for $L(w)$ if and only if $w$ does not belong to the left 
cells containing one of the following involutions: 
$s_1s_3$, $s_2s_4$, $s_2s_3s_2$, $s_1s_2s_1s_4$ or $s_1s_3s_4s_3$.
\end{proposition}

\begin{proof}
The group $S_5$ has $26$ involutions. As above, Kostant's problem
has the positive answer for $L(e)$ since it is a quotient of the 
dominant Verma module. The answers for $L(s_1)$, $L(s_2)$, $L(s_3)$ 
and $L(s_4)$ are also positive by Corollary~\ref{cor32}, and for
$L(s_1s_2s_1s_3s_2s_1s_4s_3s_2s_1)$ the answer is positive as 
this module is a Verma module. The involutions
\begin{displaymath}
\begin{array}{lll}
s_1s_2s_1,& s_1s_2s_1s_3s_2s_1, & s_1s_2s_3s_4s_3s_2s_1,\\
s_3s_4s_3,& s_2s_3s_2s_4s_3s_2, & s_2s_1s_3s_2s_1s_4s_3s_2,\\
s_3s_2s_4s_3,& s_1s_3s_2s_1s_4s_3, & s_1s_2s_3s_2s_4s_3s_2s_1,\\
s_2s_1s_3s_2,& s_2s_1s_3s_4s_3s_2, &s_1s_2s_1s_3s_4s_3s_2s_1.
\end{array}
\end{displaymath}
are all in left cells containing elements on the form $w'_0w_0$ where 
$w'_0$ is the longest element of some parabolic subgroup of $W$. Hence
Kostant's problem has the positive answer for the corresponding simple 
modules by  \cite[Theorem~4.4]{GJ} and \cite[Theorem~60]{MS1}. 
The involutions $s_2s_3s_4s_3s_2$ and $s_2s_4s_3s_2s_1$ are both in 
left cells containing elements on the form $sw'_0w_0$, where $w'_0$ 
is the longest element of some parabolic subgroup, and $s$ is a 
simple reflection of the same parabolic subgroup, so Kostant's problem 
has the positive answer for $L(s_2s_4s_3s_2s_1)$ and $L(s_2s_3s_4s_3s_2)$ 
by \cite[Theorem~1]{Ma} and \cite[Theorem~60]{MS1}.
Kostant's problem has the positive answer for $L(s_1s_4)$, and the
negative answer for $L(s_1s_3)$ and $L(s_2s_4)$, by Corollary~\ref{cor33}.

Finally, the fact that Kostant's problem has the negative answer
for $L(s_2s_3s_2)$, $L(s_1s_3s_4s_3)$ and $L(s_1s_2s_1s_4)$ follows from Theorem~\ref{thm555} by a direct computation as described in
Remark~\ref{rem557}. Consider first the involution $s_2s_3s_2$ for which
we have $\mathfrak{a}(s_2s_3s_2)=3$. A direct calculation shows that the 
graded component $\mathtt{P}^{\mathbf{\hat R}}(s_2s_3s_2)_2$
has, after forgetting the grading, the following form:
\begin{multline*}
L(s_3s_2)\oplus L(s_3s_2s_4s_3)\oplus
L(s_2s_1s_3s_2s_4s_3)\oplus  
L(s_3s_2s_1s_4s_3s_2)\oplus \\ \oplus L(s_2s_3s_2s_1)\oplus
L(s_2s_3)\oplus L(s_2s_1s_3s_2)\oplus L(s_2s_3s_2s_4).
\end{multline*}
Another calculation shows that  the graded component $\Delta(e)_2$, 
after forgetting the grading, the following form:
\begin{multline*}
L(s_3s_4)\oplus L(s_2s_4)\oplus L(s_2s_1)\oplus L(s_3s_2)\oplus
L(s_1s_3)\oplus L(s_1s_4)\oplus \\ \oplus L(s_4s_3)\oplus L(s_1s_2)\oplus
L(s_2s_3)\oplus L(s_2s_1s_3s_2)\oplus L(s_3s_2s_4s_3).
\end{multline*}
Hence the module $L(s_3s_2s_1s_4s_3s_2)$ occurs in 
$\mathtt{P}^{\mathbf{\hat R}}(s_2s_3s_2)_2$ but not in $\Delta(e)_2$. 
By Theorem~\ref{thm555} and Remark~\ref{rem556} this implies that 
Kostant's  problem has the negative answer for $L(s_2s_3s_2)$.

For the involution $s_1s_2s_1s_4$ we have $\mathfrak{a}(s_1s_2s_1s_4)=4$. 
A direct calculation shows that the module $L(s_1s_4s_3s_2s_1)$ occurs in 
$\mathtt{P}^{\mathbf{\hat R}}(s_1s_2s_1s_4)_3$ but not in $\Delta(e)_3$. Hence 
again Remark~\ref{rem556} implies that Kostant's problem has the 
negative answer for $L(s_1s_2s_1s_4)$. Applying the symmetry of the
root system we obtain that the answer for $L(s_4s_3s_4s_1)$ is
also negative and it remains to observe that $s_4s_3s_4s_1=s_1s_3s_4s_3$.
\end{proof}

Figures~\ref{fig1}, \ref{fig2} and \ref{fig3} show the three two-sided 
cells of $S_5$ which contain left cells for elements of which Kostant's 
problem has the negative answer. These left cells are columns, which are
marked with an arrow. The rows and columns are indexed by the left and 
the right Young tableaux in the corresponding Robinson-Schensted pair. 
Also, each element is denoted simply by the sequence of indices in some
shortest expression, i.e. $s_1s_3s_2$ is denoted by $132$. There seems
to exist some hidden symmetry in these pictures, but we do not understand
it yet.

\newbox\youngbox\newdimen\totheight\newdimen\raiseheight
\newdimen\ruleheight\newdimen\dropheight
\newcommand{\setpadyoungdimensions}[1]{%
\setbox\youngbox=\hbox{#1}%
\totheight=\ht\youngbox%
\raiseheight=\ht\youngbox\divide\raiseheight by 2\advance
\raiseheight by -0.7ex%
\ruleheight=\ht\youngbox\advance\ruleheight by 1ex%
\dropheight=-\raiseheight\advance\dropheight by -0.5ex%
}
\newcommand{\padyoungtab}[1]{%
\setpadyoungdimensions{#1}%
\rule[\dropheight]{0pt}{\ruleheight}\raise-\raiseheight\box
\youngbox
}

\begin{figure}
{\small\[
\begin{array}{c|c|c|c|c|c|}
& \hbox{\tiny\young(123,45)} & \hbox{\tiny\young(124,35)} &
\hbox{\tiny\young(134,25)} & \hbox{\tiny\young(135,24)} &
\hbox{\tiny\young(125,34)} \\
\hline
\padyoungtab{\tiny\young(123,45)} & 3243 & 324 & 3214 & 32143 & 321432                                             \\
\hline                              
\padyoungtab{\tiny\young(124,35)} & 243 & 24 & 214 & 2143 & 21432 \\
\hline                              
\padyoungtab{\tiny\young(134,25)} & 1321 & 124 & 14 & 143 & 1432 \\
\hline                              
\padyoungtab{\tiny\young(135,24)} & 13243 & 1324 & 134 & 13 & 132 \\
\hline                              
\padyoungtab{\tiny\young(125,34)} & 213243 & 21324 & 2134 & 213 & 2132                                            \\
\hline
\multicolumn{2}{c}{}&\multicolumn{1}{c}{\uparrow}&
\multicolumn{1}{c}{}&\multicolumn{1}{c}{\uparrow}
\end{array}
\]}
\caption{}
\label{fig1}
\end{figure}

\begin{figure}
{\small\[
\begin{array}{c|c|c|c|c|c|c|}
& \hbox{\tiny\young(123,4,5)} & \hbox{\tiny\young(124,3,5)} &
\hbox{\tiny\young(125,3,4)} & \hbox{\tiny\young(134,2,5)} & 
\hbox{\tiny\young(135,2,4)} & \hbox{\tiny\young(145,2,3)} \\
\hline
\padyoungtab{\tiny\young(123,4,5)} & 343 & 3432 & 32432 & 34321 & 
324321 & 3214321 \\
\hline                                                                 
\padyoungtab{\tiny\young(124,3,5)} & 2343 & 23432 & 2432 & 234321 & 
24321 & 214321 \\
\hline                                                                 
\padyoungtab{\tiny\young(125,3,4)} & 23243 & 2324 & 232 & 23214 & 
2321 & 21321 \\
\hline                                                                 
\padyoungtab{\tiny\young(134,2,5)} & 12343 & 123432 & 12432 & 1234321 
& 124321 & 14321                                                                \\
\hline                                                                 
\padyoungtab{\tiny\young(135,2,4)} & 123243 & 12324 & 1232 & 123214 
& 12321 & 1321 \\
\hline                               
\padyoungtab{\tiny\young(145,2,3)} & 1213243 & 121324 & 12132 & 
12134 & 1213 & 121 \\
\hline
\multicolumn{3}{c}{}&\multicolumn{1}{c}{\uparrow}
\end{array}
\]}
\caption{}
\label{fig2}
\end{figure}

\begin{figure}
{\small\[
\begin{array}{c|c|c|c|c|c|}
& \hbox{\tiny\young(12,34,5)} & \hbox{\tiny\young(12,35,4)} 
& \hbox{\tiny\young(13,24,5)} & \hbox{\tiny\young(13,25,4)} 
& \hbox{\tiny\young(14,25,3)} \\
\hline
\padyoungtab{\tiny\young(12,34,5)} & 213432 & 2132432 & 21343 
& 21324321 & 2134321                                                                       \\
\hline                                                                 
\padyoungtab{\tiny\young(12,35,4)} & 2321432 & 21321432 & 232143 
& 2132143 & 213214                                                                           \\
\hline                                                                 
\padyoungtab{\tiny\young(13,24,5)} & 13432 & 132432 & 1343 & 
1324321 & 134321 \\
\hline                                                                 
\padyoungtab{\tiny\young(13,25,3)} & 12321432 & 1321432 & 1232143 
& 132143 & 13214\\
\hline                                    
\padyoungtab{\tiny\young(14,25,3)} & 1213432 & 121432 & 121343 & 12143 
& 1214 \\
\hline
\multicolumn{3}{c}{}&\multicolumn{1}{c}{\uparrow}&\multicolumn{1}{c}{}&
\multicolumn{1}{c}{\uparrow}
\end{array}
\]}
\caption{}
\label{fig3}
\end{figure}

\subsection{Kostant's problem for $\mathfrak{sl}_6$}\label{s3.7}

We are not able yet to give a complete answer to Kostant's problem in the 
case  $\mathfrak g = \mathfrak{sl}_6$. The group $S_6$ has $76$ involutions. 
For $47$ involutions one can use arguments analogous to the
arguments above to show that Kostant's problem has the positive answer, 
for $20$ involutions one can analogously show that Kostant's problem
has the negative answer. This leaves $9$ involutions for which the answer 
is still unclear.

There are $44$ involutions, which lie in left cells containing an element 
of the form $w'_0w_0$ or $sw'_0w_0$, and hence Kostant's problem has 
the positive answer for these involutions. The remaining three for 
which Kostant's problem for sure has a positive answer are 
$L(s_1s_4)$, $L(s_1s_5)$  and $L(s_2s_5)$ (this follows from Corollary~\ref{cor33}). By the same corollary, Kostant's problem has 
the negative answer for $L(s_1s_2)$, $L(s_2s_4)$ and $L(s_3s_5)$. 
This and computations as described in Remark~\ref{rem557} show that 
Kostant's problem has the negative answer for the following $17$ involutions: 
\begin{displaymath}
\begin{array}{llll}
s_1s_3,&  s_1s_3s_5,& s_1s_4s_3s_5s_4,& s_1s_2s_1s_4s_5s_4,\\
s_3s_5,& s_1s_2s_1s_4 & s_2s_1s_3s_2s_5,& s_1s_2s_1s_3s_2s_1s_5,\\
s_2s_4,& s_1s_3s_4s_3,& s_1s_2s_3s_2s_1s_5,& s_1s_3s_4s_3s_5s_4s_3, \\
s_2s_3s_2,& s_2s_4s_5s_4,& s_1s_3s_4s_5s_4s_3,& s_1s_3s_2s_1s_4s_5s_4s_3,\\
s_3s_4s_3,& s_2s_3s_2s_5,& s_2s_3s_2s_4s_3s_2, &
s_1s_2s_1s_3s_4s_3s_5s_4s_3s_2s_1, 
\end{array}
\end{displaymath}
The remaining $9$ involutions,
which are not covered by Theorem~\ref{thm555} are:
\begin{displaymath}
\begin{array}{lll}
s_3s_2s_4s_3, & s_2s_3s_4s_3s_2, & s_2s_3s_2s_4s_5s_4s_3s_2, \\ 
s_1s_4s_5s_4, & s_2s_1s_4s_3s_2s_5s_4, & s_1s_3s_2s_4s_3s_2s_1s_5s_4s_3,\\
s_1s_2s_1s_5, & s_1s_2s_3s_2s_4s_3s_2s_1, & s_2s_1s_3s_2s_1s_4s_5s_4s_3s_2.
\end{array}
\end{displaymath}
For these involutions the answer is still unclear.

\vspace{3cm}

\noindent 
J.K.: Department of Mathematics, Uppsala University, SE-751 06, Uppsala, 
SWEDEN, e-mail: {\small \tt johank@math.uu.se}, \\
web: http://www.math.uu.se/$\tilde{\hspace{1mm}}$johank/
\vspace{0.5cm}

\noindent 
V.M.: Department of Mathematics, Uppsala University, SE-751 06, 
Uppsala, SWEDEN, e-mail: {\small \tt mazor@math.uu.se}, \\
web: http://www.math.uu.se/$\tilde{\hspace{1mm}}$mazor/ \\
and Department of Mathematics, University of Glasgow, University Gardens,
Glasgow G12 8QW, UK,\\ e-mail: {\small \tt v.mazorchuk@maths.gla.ac.uk}
\vspace{0.5cm}


\begin{thebibliography}{9999}
\bibitem[BG]{BG}
I.~Bernstein, S.~Gelfand, Tensor products of finite- and 
infinite-dimensional  representations of semisimple Lie algebras.  
{\em Compositio Math.}  {\bf 41}  (1980), no. 2, 245--285.
\bibitem[BGG]{BGG}
I.~Bernstein, I.~Gelfand, S.~Gelfand, A certain category of 
${\mathfrak g}$-modules. {\em Funkcional. Anal. i Prilozen.}  
{\bf 10}  (1976), no. 2, 1--8.
\bibitem[BB]{BB}
A.~Bj{\"o}rner, F.~Brenti, {\em Combinatorics of Coxeter groups.} Graduate 
Texts in Mathematics, {\bf 231.} Springer, New York, 2005. xiv+363 pp.
\bibitem[CPS]{CPS}
E.~Cline, B.~Parshall, L.~Scott, Finite-dimensional algebras and highest 
weight categories.  {\em J. Reine Angew. Math.}  {\bf 391}  (1988), 85--99. 
\bibitem[Di]{Di}
J.~Dixmier, Enveloping algebras. Graduate Studies in Mathematics, {\bf 11.} 
American Mathematical Society, Providence, RI, 1996. xx+379 pp.
\bibitem[GJ1]{GJ0}
O.~Gabber, A.~Joseph, Towards the Kazhdan-Lusztig conjecture.  
{\em Ann. Sci. {\'E}cole Norm. Sup.} (4)  {\bf 14}  (1981), no. 3, 261--302.
\bibitem[GJ2]{GJ}
O.~Gabber, A.~Joseph, On the Bernstein-Gelfand-Gelfand resolution and 
the Duflo  sum formula.  {\em Compositio Math.}  {\bf 43}  (1981), 
no. 1, 107--131.
\bibitem[Ir]{Ir}
R.~Irving, Projective modules in the category ${\mathcal O}_S$: 
self-duality.  {\em Trans. Amer. Math. Soc.}  {\bf 291}  (1985),  
no. 2, 701--732.
\bibitem[Ja]{Ja}
J.-C.~Jantzen, {\em Einh{\"u}llende Algebren halbeinfacher Lie-Algebren.} 
Ergebnisse der Mathematik und ihrer Grenzgebiete (3). Springer-Verlag, 
Berlin, 1983.
\bibitem[Jo]{Jo2}
A.~Joseph,  Kostant's problem, Goldie rank and the Gelfand-Kirillov 
conjecture. {\em Invent. Math.}  {\bf 56}  (1980), no. 3, 191--213.
\bibitem[KL]{KL}
D.~Kazhdan, G.~Lusztig, Representations of Coxeter groups and Hecke 
algebras.  {\em Invent. Math.}  {\bf 53}  (1979), no. 2, 165--184.
\bibitem[KM1]{KM1}
O.~Khomenko, V.~Mazorchuk, Structure of modules induced from simple 
modules with minimal annihilator.  {\em Canad. J. Math.}  {\bf 56}  
(2004),  no. 2, 293--309.
\bibitem[KM2]{KM}
O.~Khomenko, V.~Mazorchuk, On Arkhipov's and Enright's functors. 
{\em Math. Z.} {\bf 249} (2005), no. 2, 357--386. 
\bibitem[Lu]{Lu}
G.~Lusztig, Cells in affine Weyl groups.  {\em Algebraic groups and related 
topics} (Kyoto/Nagoya, 1983),  255--287, Adv. Stud. Pure Math., 6, North-Holland, Amsterdam, 1985.
\bibitem[Mm]{Mm}
M.~Marietti, Boolean elements in Kazhdan-Lusztig theory. {\em J. Algebra}
{\bf 295} (2006), no. 1, 1--26.
\bibitem[Ma1]{Ma}
V.~Mazorchuk, A twisted approach to Kostant's problem.  
{\em Glasgow Math. J.}   {\bf 47}  (2005),  no. 3, 549--561.
\bibitem[Ma2]{Ma2}
V.~Mazorchuk, Some homological properties of the category $\mathcal{O}$,
{\em Pac. J. Math.} {\bf 232} (2007), no. 2, 313--342.
\bibitem[MS1]{MS1}
V.~Mazorchuk, C.~Stroppel, Categorification of (induced) cell modules 
and the rough structure of generalized Verma modules,
preprint  arXiv:math/0702811.
\bibitem[MS2]{MS2}
V.~Mazorchuk, C.~Stroppel, Categorification of Wedderburn's basis for
$\bbC[S_n]$, preprint arXiv:0708.3949.
\bibitem[MS3]{MS3}
V.~Mazorchuk, C.~Stroppel, Projective-injective modules, Serre functors 
and symmetric algebras, preprint arXiv:math/0508119, to appear in 
{\em J. Reine Ang. Math.}
\bibitem[MiSo]{MiSo}
D.~Mili{\v c}i{\'c}, W.~Soergel, The composition series of modules 
induced from Whittaker modules. {\em Comment. Math. Helv.} {\bf 72} 
(1997), no. 4, 503--520.
\bibitem[Ri]{Ri}
C.~Ringel, The category of modules with good filtrations over a 
quasi-hereditary algebra has almost split sequences. {\em Math. Z.}
{\bf 208} (1991), no. 2, 209--223.  
\bibitem[RC]{RC}
A.~Rocha-Caridi, Splitting criteria for ${\mathfrak g}$-modules induced 
from a parabolic and the Bernstein-Gelfand-Gelfand resolution of 
a finite-dimensional, irreducible ${\mathfrak g}$-module.  
{\em Trans. Amer. Math. Soc.}  {\bf 262}  (1980), no. 2, 335--366.
\bibitem[So]{So}
W.~Soergel, Kategorie $\mathcal O$, perverse Garben und Moduln {\"u}ber 
den Koinvarianten zur Weylgruppe. {\em J. Amer. Math. Soc.} {\bf 3} 
(1990), no. 2, 421--445. 
\bibitem[St]{St}
C.~Stroppel, Category ${\mathcal O}$: gradings and translation functors.  
{\em J. Algebra}  {\bf 268}  (2003),  no. 1, 301--326.
\bibitem[Vo]{Vo}
D.~Vogan, Jr. Irreducible characters of semisimple Lie groups. II. 
The Kazhdan-Lusztig conjectures.  {\em Duke Math. J.}  {\bf 46}  
(1979), no. 4, 805--859. 
\end{thebibliography}
\end{document}